\documentclass{article}
\usepackage{amsmath}
\usepackage{amssymb}
\usepackage{amsthm}
\usepackage{mathrsfs}
\usepackage{latexsym}
\usepackage[all,cmtip]{xy}

\newcommand{\hfpbar}{H\overline{\mathbb{F}}_p} %eilenbergmac-spectra
\newcommand{\hfp}{H\mathbb{F}_{p}}

\newcommand{\s}{\mathbb{S}}%sphere spectrum
\newcommand{\Sp}[1]{\Sigma^{\infty}{#1}}

\newcommand{\op}{\mathscr{O}}%commutative operad

\newcommand{\T}{\mathscr{T}}%topological spaces
\newcommand{\Tp}{\mathscr{T}_p}
\newcommand{\F}{\mathscr{F}}%free unstable Abig alg functor
\newcommand{\Fbar}{\overline{E}}
\newcommand{\A}{\mathscr{A}} %Steenrod Algebra
\newcommand{\B}{\mathscr{B}}% E-infinity steenrod Algebra.
\newcommand{\Bbar}{\overline{\mathscr{B}}}

\newcommand{\fp}{\mathbb{F}_{p}}%field with p elts.
\newcommand{\fpbar}{\overline{\mathbb{F}}_p}
\newcommand{\p}{\mathscr{P}}%class of projectives 

\newcommand{\Map}{\mathrm{Map}}

\newcommand{\Hom}{\mathrm{Hom}}
\newcommand{\Tot}{\mathrm{Tot}}
\newcommand{\Der}{\mathrm{Der}}
\newcommand{\sAlg}{\mathrm{sAlg}}
\newcommand{\sMod}{\mathrm{sMod}}
\newcommand{\Mod}{\mathrm{Mod}}
\newcommand{\Alg}{\mathrm{Alg}}
\newcommand{\Ext}{\mathrm{Ext}}
\newcommand{\Tor}{\mathrm{Tor}}

\newcommand{\sinf}{\Sigma^{\infty}}

\newtheorem{theorem}{Theorem}[section]
\newtheorem{proposition}[theorem]{Proposition}
\newtheorem{lemma}[theorem]{Lemma}
\newtheorem{corollary}[theorem]{Corollary}

\title{A comparison of spectral sequences computing unstable homotopy groups of $p$-complete, nilpotent spaces.}
\author{Jennifer French}

\begin{document}

\maketitle

\begin{abstract}
The focus of this paper is the comparison of two unstable homotopy spectral sequences --- the unstable mod-$p$ Adams spectral sequence that computes the unstable homotopy  of a $p$-complete space, and the Goerss--Hopkins spectral sequence, which computes the unstable homotopy of the space of $E_{\infty}$ maps between $\hfpbar$-algebras.  Using an adjunction between $p$-complete nilpotent spaces and $\hfpbar$-algebras, this paper shows that the unit of this adjunction provides an isomorphism between these spectral sequences.
\end{abstract}

\section{Introduction}\label{basics}
Let $\Mod_{S}$ denote the category of spectra with model structure that of $S$-modules described in \cite{mmss}, where $S$ denotes the sphere spectrum, either as a symmetric spectrum, or as in \cite{ekmm}.  Write $\Alg_{S}$ for the category of \emph{commutative} $S$-algebras, modeling $E_{\infty}$-ring spectra.  In particular, we will fix an $E_{\infty}$ operad $\op$, and consider $\op$-algebra in spectra.  Similarly, given an $S$-algebra $R$, we may consider the categories of spectra $\Mod_R$ and $\Alg_R$.  By $\Alg_R$ we will always mean \emph{commutative} $R$-algebras, the category of $S$-algebras under $R$.  Let $\hfpbar$ denote the Eilenberg--Maclane spectrum whose homotopy $\pi_*\hfpbar = \fpbar$ for $*=0$ and is 0 otherwise.  In particular, we are interested in the category $\Alg_{\hfpbar}$ of commutative $\hfpbar$-algebras, where $p$ is a prime.

In the unstable world, we are interested in $p$-complete spaces.  There are two common definitions for the $p$-completion of a space.  One being the $\hfp$ localization of the space, and the other being the Bousfield--Kan $p$-completion.  In general, these do not agree.  However, for a \emph{nilpotent} space $X$, both notions of $p$-completion coincide \cite[Proposition V.4.2]{holimbk}.  A space $X$ is \emph{nilpotent} if $\pi_1X$ is nilpotent and acts nilpotently on the higher homotopy groups $\pi_nX$ for $n>1$.  We restrict our attention to nilpotent spaces throughout this exposition.  

The \emph{$p$-completion} of a space $X$ is the homotopy inverse limit of the diagram of spaces:
\[ \xymatrix{\fp(X) \ar@<1ex>[r] \ar@<-1ex>[r] &\ar[l] \fp(\fp(X)) \ar[r]  \ar@<2ex>[r] \ar@<-2ex>[r]  &\ar@<-1ex>[l] \ar@<1ex>[l]\dotsb.}\] The maps in this cosimplicial diagram are discussed in Section \ref{uass}.

Let $\T$ denote the category of (unpointed) simplicial sets.  Note that all constructions carried out in this paper have analogs when considering $\T$ to be the category of weak Hausdorff topological spaces.  We take the model of simplicial sets in order to make certain constructions more transparent.  The author uses the convention of referring to simplicial sets as spaces.  Recall that there is a functor $(-)_+$ from unpointed to pointed spaces given by adding a disjoint basepoint, which is left adjoint to the forgetful functor.  Let $\Tp$ denote the category of $p$-complete, nilpotent spaces.

Note that $\Alg_{\hfpbar}$ is enriched over spaces.  And the model category structure satisfies Quillen's corner axiom.  Thus given any two $\hfpbar$-algebras $A$ and $B$, referring to the derived mapping space $\Alg_{\hfpbar}(A,B)$ makes sense.  In this paper we follow the convention that given any category $\mathscr{C}$ enriched in spaces, $\mathscr{C}(A,B)$ denotes the derived mapping space, while $\pi_0\mathscr{C}(A,B)$ denotes the underlying set of morphisms in $\mathscr{C}$.  

Consider the functor from $\T$ to $\Alg_{\hfpbar}$ discussed in \cite[Appendix C]{mandell}:
\begin{align*}
F:\T&\rightarrow \Alg_{\hfpbar},\\
X &\mapsto \hfpbar^{\Sigma^{\infty}X_+}.
\end{align*} 
This is a functor of simplicial categories.  In \cite{mandell}, M. Mandell proves that $F$ induces an equivalence onto its image when restricted to the subcategory of $p$-complete, nilpotent spaces.  In particular, the homotopy inverse of $F$ is given by the right adjoint \begin{align*} G: \Alg_{\hfpbar} &\rightarrow \T,\\ A &\mapsto \Alg_{\hfpbar}(A, \hfpbar). \end{align*}  Thus, given a nilpotent space, the composite $\psi := G \circ F : \T \rightarrow \T$ induces an $\hfp$-localization:
\begin{align*}
\psi :  X \rightarrow \Alg_{\hfpbar}(\hfpbar^{\Sp{X_+}}, \hfpbar).
\end{align*} 

More is true.  Applying the functor $\hfpbar^{\sinf (-)_+}$ to an unstable map of simplicial sets gives rise to a map of $\hfpbar$-algebras.  Thus given a mapping space $\T(Y,X)$, there is an induced functor to the mapping space of $\hfpbar$-algebras.  The theorem as proved by M. Mandell in \cite{mandell} can be rewritten in the language of $\hfpbar$-algebras as the following:

\begin{theorem}[\cite{mandell}]\label{mandell} Let $X$ and $Y$ be nilpotent and of finite type.  Then the map 
\begin{align}\label{mainmap}
\psi: \T(Y,X_{\hfp})\rightarrow \Alg_{\hfpbar}(\hfpbar^{\sinf X_+}, \hfpbar^{\sinf Y_+})
\end{align} is a weak equivalence.\end{theorem}

The unstable homotopy groups of a ($p$-complete) space $X$ can be computed using the unstable Adams spectral sequence.  The homotopy groups of the space of $\hfpbar$-algebra maps can be computed using a Goerss--Hopkins spectral sequence.  The main question we set out to answer is whether or not $\psi$ induces an isomorphism between these two spectral sequences.  The answer is yes and is the main application of the main theorem in this paper.

The goal of this exposition is to prove a general theorem and arrive at the isomorphism of the two unstable spectral sequences as a corollary.  There is a canonical simplicial resolution $P(F(\sinf X_+, \hfpbar))_{\bullet}$ given by the cofibrant replacement of the constant simplicial object in the $E_2$ model category structure on simplicial $\hfpbar$-algebras.  The cosimplicial space $\fp^{\bullet}(X)$ is obtained by iteratively applying a the monad $\fp$ sending $X$ to the underlying simplicial set on the free simplicial $\fp$-algebra on $X$.  Let $Y$ be a space and consider the diagram of cosimplicial spaces:
\begin{align*}
\xymatrix{ \T(Y,\fp^{\bullet}X) \ar[r] & \mathrm{diag}\Alg_{\hfpbar}(P(F(\sinf \fp^{\bullet}X_+,\hfpbar^{\sinf Y_+}))_{\bullet},\hfpbar) \\ & \ar[u] \Alg_{\hfpbar}(P(F(\sinf X_+, \hfpbar))_{\bullet},\hfpbar^{\sinf Y_+}).}\end{align*}  The main theorem is that both maps are equivalences when $X$ and $Y$ are of finite type.  In particular, this theorem offers another proof of Theorem \ref{mandell} upon totalization of this cosimplicial diagram.  It should be pointed out that although this proof is different, the main computation used in proving the main theorem is the same as in M. Mandell's original proof.  The corollary exhibits an isomorphism between the $E_2$ terms of the unstable $\hfp$-Adams spectral sequence and the Goerss--Hopkins spectral sequence when conditions of the theorem are satisfied.

In Section \ref{unstablify} we define two categories of unstable algebras, over the Steenrod algebra and the $E_{\infty}$ Steenrod algebra.  We define Andr\'e--Quillen cohomology theories in these categories.  The purpose of Section \ref{descent} is a natural isomorphism between the derived functors of derivations of unstable $\fpbar$-algebras over the $E_{\infty}$ Steenrod algebra with the derived functors of derivations of unstable algebras over the Steenrod algebra when the algebras in question have a particular form.  This isomorphism is the main tool we use to identify the unstable Adams spectral sequence with the Goerss--Hopkins spectral sequence from the $E_2$-term onward.  In Section \ref{ua} we recall the unstable Adams spectral sequence and the identification of its $E_2$ term as an Andr\'e--Quillen cohomology group of unstable algebras over the mod-$p$ Steenrod algebra.  In Section \ref{gh} we recall the Goerss--Hopkins spectral sequence and identify its $E_2$ term algebraically as the Andr\'e--Quillen cohomology of unstable algebras over the $E_{\infty}$ Steenrod algebra.  In Section \ref{alg} we begin to analyze a map between two spectral sequences induced by applying the map $\psi$ from Theorem \ref{mandell} to the standard cosimplicial resolution of a space $X$ constructed in Section \ref{ua}.  The work in this section proves that $\pi^s\pi_t$ of the image under application of the map $\psi$ to the Bousfield--Kan cosimplicial $\fp$ resolution of a space $X$ can in fact be identified as the $E_2$ term of the Goerss--Hopkins spectral sequence.    The cosimplicial space we obtain via applying $\psi$ to this canonical resolution is the derived mapping space out of a simplicial $\hfpbar$-algebra.  This simplicial $\hfpbar$-algebra can be expressed levelwise as the homotopy cofiber of free algebras over the commutative operad.  Forming a bisimplicial object interpolating between the simplicial object we obtain via $\psi$ and a cofibrant object in the model structure used to construct the Goerss--Hopkins spectral sequence, we obtain a collapsing spectral sequence, which gives a levelwise weak equivalence of cosimplicial spaces.  This allows us to identify the image of $\psi$ with the $E_2$ term of the Goerss--Hopkins spectral sequence, and show that $\psi$ induces an isomorphism of spectral sequences from the $E_2$ term onward.  

I would like to thank my advisor Mark Behrens.  I am indebted to his many helpful conversations, guidance through the literature, and careful editing of this document.  I would also like to thank the referee for many helpful comments and suggestions for improving this paper.

\section{Unstable algebras over $\A$ and $\B$}\label{unstablify}
The two spectral sequences we wish to compare have $E_2$ terms that can be identified as certain Andr\'e--Quillen cohomology groups.  When restricted to spaces of finite type, the $E_2$-term of the  unstable Adams spectral sequence can be identified as derived functors of derivations of unstable algebras over the Steenrod algebra $\A$.  The $E_2$-term of the Goerss--Hopkins spectral sequence is naturally identified as derived functors of derivations of unstable algebras over the $E_{\infty}$ Steenrod algebra $\B$ over $\fpbar$.  

The goal of this section is to describe these categories of unstable algebras and Andr\'e--Quillen cohomology in these categories.  This sets up the necessary background for Section \ref{descent} where we give a descent condition relating the Andr\'e--Quillen cohomology theories over these categories.

\subsection{The Steenrod algebra}
Let $\A$ denote the Steenrod algebra over $\fp$.  We recover $\A$ as the ring $\hfp^*\hfp$ of $\hfp$ cohomology operations.  Recall that at an odd prime $p$, the ring $\A$ over $\fp$ is the quotient of the free $\fp$ algebra generated by $P^i$ and $\beta P^i$, for $i$ any non-negative integer, modulo the ideal generated by the Adem relations.  At the prime 2 we write $\A$ as the quotient of $\mathbb{F}_{2}[Sq^0, Sq^1, \dotsc]$ by the ideal generated by the Adem relations \cite{schwartz}.

There is a basis for $\A$ given by \emph{admissible} monomials.  Let the symbols $\epsilon_k \in \{0,1\}$.  A sequence $I =  (\epsilon_0, i_0, \dotsc ,\epsilon_n, i_n)$ is admissible if $i_{k-1} \leq pi_k +\epsilon_k$.  The $\fp$-module basis for $\A$ is given by $\{P^I: I \mathrm{admissible}\}$, where $P^{\epsilon_i}$ denotes $\beta^{\epsilon_i}$.  

An $\A$-module is a $\fp$ vector space $M$ with an action of $\A$.  Define the excess of a monomial $I$ to be  \begin{align*} e(I)= & i_1+ - \sum_{k=2}^{n}i_k & \qquad p = 2\\
 e(I)= & 2i_1+\epsilon_1 - \sum_{k=2}^{n}(2(p-1)i_k +\epsilon_k)& \qquad p \neq 2.\end{align*}  Then $M$ is an \emph{unstable} $\A$-module if for all $x$ in $M$, $P^I(x) = 0$ whenever $e(I)>|x|$.  There is a Quillen pair of adjoint functors from graded $\fp$ vector spaces $gr\Mod_{\fp}$ to unstable $\A$ modules $\Mod_{U\A}$
\begin{align*}\xymatrix{
\F_0: gr\Mod_{\fp} \ar@<.5ex>[r] & \ar@<.5ex>[l] \Mod_{U\A}: U
}\end{align*}
where $U$ simply forgets the unstable $\A$-module structure, and the left adjoint $\F_0$ is the free functor.  We can define $\F_0$ explicitly by choosing an $\fp$ basis $\{x_j\}_{j\in J}$ of our vector space $M$.  Then $\F_0(M)$ is the unstable $\A$ module with $\fp$-module basis given by $\{P^Ix_j : j \in J, I \mathrm{ admissible}\}.$  The main example of an unstable $\A$-module appearing in nature is the homotopy groups of the $\hfp$-module $F(A,\hfp)$, where $A$ is a -1-connected spectrum.

An unstable $\A$-module $B$ is an unstable $\A$-algebra if is a ring, and satisfies the instability condition $P^ix = x^p$ whenever $i = |x|$.  That is, $P^i$ acts as the frobenius on the $i$th graded piece of $B$.  The natural example of an unstable $\A$-algebras is $\hfp$ cohomology ring of a space $X$. 

There is a Quillen pair of adjoint functors from unstable $\A$-modules to unstable $\A$-algebras:
\begin{align*}
\xymatrix{E_0: \Mod_{U\A} \ar@<.5ex>[r] & \ar@<.5ex>[l] \Alg_{U\A}: U,}
\end{align*} where the right adjoint $U$ is the forgetful functor, and the left adjoint $E_0$ takes the free commutative algebra on a module $M$ and identifies $P^i m$ with $m^p$ whenever the degree of $m$ is $i$.

Use the composition of left Quillen functors to define a model category structure on the category of simplicial unstable $\A$-algebras so that the cofibrant objects are given by simplicial unstable $\A$-algebras that are levelwise free.

\begin{lemma}
Suppose we have unstable $\A$-algebras $A_*$ and $B_*$, a map $\varphi : A_* \rightarrow B_*$, and an unstable $\A$-module, $M_*$, that is also $B_*$-module.  If $\A$ acts trivially on $M_*$, then $B_* \ltimes M_*$ is naturally an unstable $\A$-algebra over $B_*$.  There is an isomorphism \[\Alg_{U\A \downarrow B_ *}(A_*, B_*\ltimes M_*) \cong \Der_{U\A/B_*}(A_*,M_*),\] where $\Der_{U\A/B_*}$ is defined to be $\fp$- linear derivations relative to $B_*$ that commute with $\A$. 
\end{lemma}
\begin{proof}
First, $B_*\ltimes M_*$ is an unstable $\A$ algebra with action $P^i(b,m) = (P^ib, P^im) = (P^ib, 0)$.  Multiplication is given by \begin{align*}
(b,m)\otimes (b', m') & \mapsto (bb', b\otimes m'+b' \otimes m).\end{align*} 

A map $f$ in $\Alg_{U\A \downarrow B_*}(A_*, B_* \ltimes M_*)$ is a pair of maps $f = (\varphi, g)$ such that \[ f(aa') = (\varphi(aa'), \varphi(a)g(a') + \varphi(a')g(a)),\] commuting with the action of $\A$.  In particular, $P^if = (P^i\varphi, g)$.  The map $g: A_* \rightarrow M_*$ is a derivation with trivial $P^i$.  This follows from the definition of the multiplication in $B_* \ltimes M_*$
\begin{align*}
g(aa') & = \varphi(a)g(a') + \varphi(a')g(a),\\
\end{align*} and the fact that $\A$ acts trivially on $M_*$.  

Now, given a derivation $d: A_* \rightarrow M_*$ commuting with $\A$, it is clear that the map
\[ \xymatrix { 
A_* \ar[rr]^{(\varphi,d)} \ar[dr]_{\varphi} && B_* \ltimes M_* \ar[dl]^{p_1}\\
 &B_*
}\] is a map of unstable $\B$-algebras over $B_*$:
\begin{align*}
ab & \mapsto (\varphi(ab), \varphi(a)d(b) + \varphi(b)d(a))\\
P^i(ab) & \mapsto (\varphi(P^i(ab)), \Sigma_{i_0+i_1=i}(P^{i_0}\varphi(a)P^{i_1}d(b)+P^{i_0}\varphi(b)P^{i_1}d(a)) \\
& \qquad \cong (\varphi(P^i(ab)),\varphi(P^ia)d(b)+\varphi(P^ib)d(a)).  
\end{align*}
\end{proof}

Thus the derived functors of derivations are computed by taking a cofibrant replacement for $A_*$ in the category of simplicial unstable $\A$-algebras.  Such a cofibrant replacement is a free unstable $\A$-algebra resolution $F(A_*)_{\bullet} \rightarrow A_*$.  Define the Andr\'e--Quillen cohomology as:
\begin{align}
D^*_{U\A/B_*}(A_*,M_*) := \pi^*\Der_{U\A/B_*}(F(A_*)_{\bullet},M_*).
\end{align}

\subsection{The $E_{\infty}$ Steenrod algebra}\label{steenrodalg}
Let $\B$ denote the $E_{\infty}$ Steenrod algebra over $\fp$.  Recall that at an odd prime $p$, the ring $\B$ over $\fp$ is the quotient of the free $\fp$ algebra generated by $P^i$ and $\beta P^i$ for $i$ any integer modulo the ideal generated by the Adem relations.  At the prime 2 we write $\B$ as the quotient of $\mathbb{F}_{2}[\dotsc,Sq^{-1}, Sq^0, Sq^1, \dotsc]$ modulo the ideal generated by the Adem relations \cite{mandell}.  

A sequence $I = (\epsilon_1, s_1, \epsilon_2, s_2,\dotsc,\epsilon_k, s_k)$ is called admissible if $s_{j-1} \leq ps_j +\epsilon_j.$  There is an additive basis for $\B$ given by $\{P^I: I = (i_0, \epsilon_0, i_1, \epsilon_1,\dotsb) \mathrm{ admissible}\}$.  Recall that $\B$ is a graded ring with $|P^i| = 2i(p-1)$ and B\"ockstein $|\beta|=1$.  
    
Given an admissible sequence $I$, let \begin{align*}e(I)=  & s_1+ - \sum_{j=2}^{k}s_j & \qquad p=2\\
 e(I)= & 2s_1+\epsilon_1 - \sum_{j=2}^{k}(2(p-1)s_j +\epsilon_j) & \qquad p\neq2 \end{align*} denote the excess of $I$.  We say that a module $M$ over $\B$ is unstable if for all $x$ in $M$, the formula $P^I(x) =0$ whenever $e(I)> |x|$.  There is a Quillen pair of adjoint functors from graded $\fp$ vector spaces $gr\Mod_{\fp}$ to unstable $\B$ modules $\Mod_{U\B}$
\begin{align*}\xymatrix{
\F: gr\Mod_{\fp} \ar@<.5ex>[r] & \ar@<.5ex>[l] \Mod_{U\B}: U
}\end{align*}
where $U$ simply forgets the unstable $\A$-module structure, and the left adjoint $\F_0$ is the free functor.  We can define $\F$ explicitly by choosing an $\fp$ basis $\{x_j\}_{j\in J}$ of our vector space $M$.  Then $\F(M)$ is the unstable $\B$ module with $\fp$-module basis given by $\{P^Ix_j : j \in J, I \mathrm{ admissible}\}.$  The main example of an unstable $\B$-module appearing in nature is the homotopy groups of a (non-connective) $\hfp$-module. 

An unstable module is an unstable algebra if the algebra map is a map of unstable $\B$-modules, and $P^ix = x^p$ whenever $i= |x|.$

Let $A$ be an unstable $\B$-algebra.  An unstable $\B$ module $M$ is an $A$-module if there is $\B$-module map $A \otimes M \rightarrow M$ so that the induced action satisfies the Cartan relation 
\begin{align*}
P^i(a m) = \Sigma_{i_0 + i_1 = i} P^{i_0}a P^{i_1}m
\end{align*} where $a \in A$ and $m \in M$.  Observe that the instability condition assures that this sum is finite.

Every unstable $\B$-module is naturally a graded restricted $\fp$-module --- a graded $\fp$-module with a restriction map $\Phi: M_{*} \rightarrow M_{p*}$ that multiplies degree by $p$.  In the case that the restricted module comes from an unstable $\B$-module, at odd primes $p$, the restriction map is given by $P^m: M_m \rightarrow M_{pm}$ in even degrees, and is identically zero on odd degrees.  
There is a Quillen pair of functors from unstable $\B$-modules to unstable $\B$-algebras
 \begin{align}\xymatrix{
 E: \Mod_{U\B} \ar@<.5ex>[r] &  \ar@<.5ex>[l] \Alg_{U\B}:U.
 }\end{align}  The left adjoint $E$ is given by the enveloping algebra functor on the underlying graded restricted $\fp$-module.  This has the effect of equating the restriction map with the frobenius.  The right adjoint is the forgetful functor.
 
\subsection{Relations between $\A$ and $\B$} 
\begin{lemma}[{\cite[Proposition 12.5]{mandell}}] \label{sesresmod}
For any graded $\fp$ vector space $V$, there is a natural short exact sequence of $\B$-modules:
\[ \xymatrix{ 0 \ar[r] & \F(V) \ar[r]^{1-P^0} & \F(V) \ar[r]^q & \F_0(V) \ar[r] & 0,} \] which is split exact on the underlying graded restricted $\fp$-modules.
\end{lemma}
Note that in particular, this Lemma says that we can recover the Steenrod algebra $\A$ from $\B$ by taking the quotient of $\B$ by the left ideal (which is incidentally a two-sided ideal) generated by $(1-P^0)$ {\cite[Proposition 11.4]{mandell}}.   Thus we can view an unstable $\A$-module as an unstable $\B$-module via the quotient map.  

Recall that the homotopy groups of any $\hfp$-algebra naturally form an unstable $\B$-algebra.  But in the case that the $\hfp$-algebra is free over $\op$, the homotopy groups naturally form a free unstable $\B$-algebra.

\begin{lemma}
If $Z$ is an $\hfp$-module, then \[ \pi_*\op(Z) \cong E\F(\pi_*Z).\]
\end{lemma}
The proof in the case of the Steenrod algebra can be found in \cite{schwartz} or in the case of the Dyer-Lashoff algebra in \cite[Theorem IX.2.1]{h-infty}.  The proof for the $E_{\infty}$ Steenrod algebra follows from these proofs.

Applying the previous two lemmas, we find that $\pi_*\op \hfp^{S^n_+} \cong E\F(\fp[-n]\oplus \fp[0]).$  Since we know that $\hfp^*K(\fp,n)$ is the free unstable $\A$ algebra on a generator of cohomological degree $n$, the we find the following natural extension of Lemma \ref{sesresmod}.
\begin{lemma}\label{spushout}
The map \[ 1-P^0: E\F(\fp[-n]) \rightarrow E\F(\fp[-n])\] is injective and the target is a projective module over the source via this map.  Moreover, there is a pushout diagram in the category of unstable $\B$-algebras:
\[ \xymatrix{  E\F(\fp[-n]) \ar[d] \ar[r]^{1-P^0} & E\F(\fp[-n])  \ar[d] \\
\fp \ar[r] & E\F_0(\fp[-n]).
}\] 
\end{lemma}
Proof of this Lemma can be found in \cite{mandell}.  It is important to note that $E\F_0(V) \cong E_0\F_0(V)$ since by definition the Dyer--Lashoff operations in $\B$ act trivially on $\F_0(V).$

In fact, even more is true.  Let $\fp[-n]$ denote the graded $\fp-$module with one $\fp$ in cohomological degree $n$.  Then we use notation from \cite{mandell} and denote $\F(\fp[-n])$, the free $\B-$module on a generator in degree $n$, as $\B^n$.  As an example, the enveloping algebra on a free unstable $\B-$module $\B^n$ is given by the free commutative algebra on the symbols $\{P^Ii_n: e(I) < n\}$ and $P^{n}i_n = i_n^p,$ where $i_n$ denotes the generator of $\B^n$ as a free $\B-$module.

We have a left Quillen functor from graded $\fp$-modules to unstable $\B$-modules.  By  {\cite[Theorem II.4.4]{quillen}}, there is a model category structure on simplicial unstable $\B$ modules such that the fibrations and weak equivalences are the maps that are fibrations and weak equivalences on the underlying simplicial graded $\fp$ modules. Consider the pushout diagram in chain complexes of unstable $\B$-modules:
\[\xymatrix{
\B^n[0]\ar[d] \ar[r]^{1-P^0}& \B^n[0] \ar[d]\\
0 \ar[r] & C_*,
}\] where $C_*$ is the chain complex of unstable $\B$-modules given by \[C_*= \dotsc 0 \rightarrow \B^n[1]\rightarrow \B^n[0].\]  Considering the image under the Dold--Kan correspondence, we obtain the simplicial unstable $\B$-module $T_{\bullet}$ as the homotopy pushout of simplicial unstable $\B$-modules:
\begin{align}\label{push}\xymatrix{
\B^n_{\bullet}\ar[d] \ar[r] &  \B^n_{\bullet}\ar[d]\\
0 \ar[r] & T_{\bullet}.
}\end{align}
 
\begin{proposition}\label{freeres}
There is a cofibrant simplicial resolution of $\pi_*\hfp^{K(\fp,n)_+}$ by free unstable $\B$-algebras given by applying $E$ homotopy pushout $T_{\bullet}$ of the diagram (\ref{push}):
\[ E T_{\bullet} \rightarrow \pi_*\hfp^{K(\fp,n)_+}.\]
\end{proposition}
\begin{proof}
There is an augmentation map $T_{\bullet} \rightarrow \F(\fp[-n]),$ which is a weak equivalence of simplicial $\B$-algebras by viewing $\F(\fp[-n])$ as a constant simplicial object.  Then $E(T_{\bullet})$ is the bar construction on the map: \[ E(1-P^0): E\F(\fp[-n]) \rightarrow  E\F(\fp[-n]).\]  Thus the homotopy groups of $ET_{\bullet}$ are given by \[ \pi_*ET_{\bullet} = \Tor_{*}^{E\F(\fp[-n])}(\fp,E\F(\fp[-n])).\]  Since $E\F(\fp[-n])$ is projective as a $E\F(\fp[-n])$ module via the module map $E(1-P^0)$ by Proposition \ref{spushout}, the higher homotopy groups vanish.  Thus \[\pi_*ET_{\bullet} \cong E \B^n \otimes_{E \B^n} \fp.\]  But we already know from \cite{mandell} that this right hand side is $\hfp^*K(\fp,n),$ which is what we wanted to show.
\end{proof}

\subsection{Unstable $\Bbar$-algebras.}
Now we extend the discussion above to include $\hfpbar$-algebras and unstable $\B$-algebras over $\fpbar$, which we will denote as $\Bbar$-algebras.  This is the context  which is important for identifying the $E_2$ term of the Goerss--Hopkins spectral sequence.

Let $f:\fpbar \rightarrow \fpbar$ denote the frobenius map on $\fpbar$.  Note that this map is injective and$f$ is the generator of the Galois group of $\fpbar$ over $\fp$.  An unstable $\B$-algebra over $\fpbar$ is an unstable $\B$ algebra $A$, whose underlying graded vector space is an $\fpbar$-module.  Given any $a \in A$ and any $\lambda \in \fpbar$, the action of an element $P^i$ is defined as:
\[P^i(\lambda a) := f(\lambda)P^ia.\]  In particular, $P^0$ acts as the frobenius on $\fpbar$.  An example of unstable $\B$-algebras over $\fpbar$ is the homotopy groups of an $\hfpbar$-algebra.  We denote the category of unstable $\Bbar$-algebras as $\Alg_{U\Bbar}$.

\begin{lemma}\label{adaption}
If $Z$ is an $\hfpbar$-module, then \[ \pi_*\op(Z) \cong E(\F(\pi_*Z)),\] is a free unstable $\B$-algebra over $\fpbar$.
\end{lemma}
We have the following diagram of Quillen adjunctions, where only the left adjoints are labeled:
\begin{align*}\xymatrix{
\Alg_{U\B} \ar@<.5ex>[r]^{-\hat{\otimes}\fpbar} \ar@<.5ex>[d] & \ar@<.5ex>[l]  \ar@<.5ex>[d] \Alg_{U\Bbar} \\
\Mod_{U\B} \ar@<.5ex>[r]^{-\hat{\otimes}\fpbar} \ar@<.5ex>[d]\ar@<.5ex>[u]^{E} & \ar@<.5ex>[l]  \ar@<.5ex>[d]\ar@<.5ex>[u]^{E} \Mod_{U\Bbar}\\
gr\Mod_{\fp} \ar@<.5ex>[r]^{-\hat{\otimes}\fpbar} \ar@<.5ex>[u]^{\F} & \ar@<.5ex>[l]  \ar@<.5ex>[u]^{\F} gr\Mod_{\fpbar}
}\end{align*}  Write \[\Fbar:= E(-)\hat{\otimes}_{\fp}\fpbar: \Mod_{U\B} \rightarrow \Alg_{U\Bbar}\] for the free functor from unstable $\B$-modules to unstable $\B$-algebras over $\fpbar$.  This is a left Quillen adjoint to the forgetful functor.  

Lemma \ref{adaption} above shows that the operad $\op$ in $\Mod_{\hfpbar}$ is \emph{adapted} to $\pi_*$ as in {\cite[Definition 1.4.13]{ghlong}}.  By {\cite[Example 1.4.14]{ghlong}} this means that $\op$ is simplicially adapted to $\pi_*$.  Thus, given $Z_{\bullet}$, a cofibrant simplicial $\hfpbar$-module, there is an isomorphism of simplicial unstable $\B$-algebras over $\fpbar$ \begin{align}\label{simpleres} \pi_*\op(Z_{\bullet}) \cong E\F(\pi_*Z_{\bullet}).\end{align}

We record here an extension of Lemma \ref{spushout} to finite type complexes extended over $\fpbar$.  This Proposition will be crucial in constructing the map of spectral sequences.

\begin{proposition}\label{coeq}
For every graded $\fp$ vector space of finite type, \[ \Fbar(1-P^0): \Fbar(\F(V))\rightarrow  \Fbar(\F(V))\] is injective, and $ \Fbar(\F(V))$ is projective as a module over itself via $ \Fbar(1-P^0)$.  Moreover, there is a pushout diagram in the category of unstable $\B$-algebras over $\fpbar$.
\[\xymatrix{\Fbar(\F(V)) \ar[r]^{1-P^0}\ar[d] & \Fbar(\F(V)) \ar[d] \\
\fpbar \ar[r] & \Fbar\F_0(V) }\] 
\end{proposition}

\begin{proof}  The proof follows from Lemma \ref{spushout} and Lemma \ref{adaption}.
\end{proof}

We also extend Proposition \ref{freeres} to spaces of finite type.  This resolution will be key in order to compute the $E_2$ term of the Goerss--Hopkins spectral sequence. 
\begin{proposition}\label{mainfreeres}  Let $V$ be an $\fp$ vector space.  Then viewing $\F(V)_{\bullet}$ as a constant simplicial object, let $T_{\bullet}$ denote the homotopy cofiber
\[\xymatrix{ \F(V)_{\bullet}\ar[r]^{1-P^0}& \F(V)_{\bullet} \ar[r] & T_{\bullet}.
}\]  Then 
\[\xymatrix{\Fbar(T_{\bullet}) \ar[r] & \Fbar\F_0(V)
}\] is a resolution of $\hfpbar^*\F_0(W)$ by a simplicial free $\B$-algebra over $\fpbar$. 
\end{proposition} 
\begin{proof}
The proof is immediate as $\Fbar$ is a left Quillen functor, and is completely analogous to the proof of Proposition \ref{freeres}.
\end{proof}

\subsection{Derivations of unstable $\Bbar$-algebras.}
\begin{lemma}
Given unstable $\Bbar$-algebras $A_*$ and $B_*$, a map $\varphi : A_* \rightarrow B_*$, and an unstable $\Bbar$-module, $M_*$, that is also $B_*$-module.  If the cohomological Steenrod operations act trivially on $M_*$, then $B_* \ltimes M_*$ is naturally an unstable $\Bbar$-algebra over $B_*$.  There is an isomorphism \[\Alg_{U\Bbar \downarrow N_ *}(A_*, B_*\ltimes M_*) \cong \Der_{U\Bbar/B_*}(A_*,M_*),\] where $\Der_{U\Bbar}$ is defined to be $\fpbar$- linear derivations relative to $B_*$ that commute with $\B$. 
\end{lemma}
\begin{proof}
A map $f$ in $\Alg_{U\Bbar \downarrow B_*}(A_*, B_* \ltimes M_*)$ is a pair of maps $f = (\varphi, g)$ such that \[ f(aa') = (\varphi(aa'), \varphi(a)g(a') + \varphi(a')g(a)).\]  Moreover, as this is a map of unstable $\B$-algebras, it commutes with the action of $\B$.  The map $g: A_* \rightarrow M_*$ is a derivation commuting with the $P^i$.  This follows from the string of equalities:
\begin{align*}
g(P^i(aa')) & = \Sigma_{i_0+i_1=i}g(P^{i_0}aP^{i_1}a')\\
&= \Sigma_{i_0+i_1=i} \varphi(P^{i_0}a)g(P^{i_1}a')+\varphi(P^{i_0}a')g(P^{i_1}a)\\
& = P^i(\varphi(a)g(a')+\varphi(a')g(a))\\
& = P^i(g(aa')).
\end{align*}  Now, given a derivation $d: A_* \rightarrow M_*$ commuting with $\B$, it is clear that the map
\[ \xymatrix { 
A_* \ar[rr]^{(\varphi,d)} \ar[dr]_{\varphi} && B_* \ltimes M_* \ar[dl]^{p_1}\\
 &B_*
}\] is a map of unstable $\Bbar$-algebras over $B_*$.
\end{proof}

The derived functors of derivations are computed by taking a cofibrant replacement for $A_*$ in the category of simplicial unstable $\B$-algebras over $\fpbar$.  Such a cofibrant replacement is a free unstable $\B$-algebra resolution $F(A_*)_{\bullet} \rightarrow A_*$.  Define the Andr\'e--Quillen cohomology as:
\begin{align}
D^*_{U\Bbar/B_*}(A_*,M_*) := \pi^*\Der_{U\Bbar/B_*}(F(A_*)_{\bullet},M_*).
\end{align}

\section{A descent theorem.}\label{descent}
\begin{proposition}
Let $V$ be a graded $\fp$-vector space.  Let $M$ be an $\Fbar\F_0(V)$-module in the category of unstable $\B$-algebras over $\fpbar$.  Let $\Fbar\F(V)$ be augmented over the $U\Bbar$-algebra $A$.  Then there is a long exact sequence in Andr\'e--Quillen cohomology: 
\begin{align*}
\dotsb & \rightarrow D^{s-1}_{U\Bbar/A}(\Fbar\F(V),M) \rightarrow  D^{s}_{U\Bbar/A}(\Fbar\F_0(V),M) \rightarrow \\
& D^{s-1}_{U\Bbar/A}(\Fbar\F(V),M) \rightarrow  D^{s-1}_{U\Bbar/A}(\Fbar\F(V),M) \rightarrow \dotsb.
\end{align*}
\end{proposition}
\begin{proof}
Proposition \ref{mainfreeres} gives a homotopy cofiber sequence of simplicial unstable $\Bbar$-algebras:
\begin{align}
\Fbar\F(V) \rightarrow \Fbar\F(V) \rightarrow \Fbar T_{\bullet}
\end{align} where the first two simplicial objects are constant.  Combining this with Proposition \ref{coeq}, which says that $\Fbar T_{\bullet}$ is a cofibrant simplicial resolution of $\Fbar \F_0(V)$ and the fact that the two free constant simplicial objects are already cofibrant, the result is immediate.
\end{proof}

Of course, since $\Fbar\F(V)$ is already cofibrant, this implies that the higher cohomology vanishes.  Thus this long exact sequences reduces to the exact sequence
\begin{align}\label{4term}
0& \rightarrow D^{0}_{U\Bbar/A}(\Fbar\F_0(V),M) \rightarrow  D^{0}_{U\Bbar/A}(\Fbar\F_0(V),M) \rightarrow \\
& D^{0}_{U\Bbar/A}(\Fbar\F(V),M) \rightarrow  D^{1}_{U\Bbar/A}(\Fbar\F_0(V),M) \rightarrow 0.
\end{align}

\begin{theorem}\label{cor}  Let $V_0$ be a graded $\fp$-module, and suppose $E\F_0(V_0)$ is augmented over the unstable $\A$-algebra $A_0$.  Write $A:= \fpbar \otimes_{\fp} A_0$.  Let $M_0$ be an $\Fbar\F_0(V_0)$-module in the category of unstable $\A$-algebras.  Define $V= \fpbar \otimes_{\fp}V_0$ and $M = \fpbar \otimes_{\fp} M_0.$   
\begin{enumerate}
\item $D^{1}_{U\Bbar/E\F_0(V)}(E\F_0(V),M)=0$
\item There is a natural isomorphism
\begin{align*}
D^{0}_{U\Bbar/A}(E\F_0(V),M) \cong D^{0}_{U\A/A_0}(E\F_0(V_0),M_0).
\end{align*}
\end{enumerate}
\end{theorem}

\begin{proof}
Note that $\Fbar\F(V_0) \cong E\F(V)$, which is augmented over the unstable $\Bbar$-algebra $A$.
We compute $D^0_{U\Bbar/A}(E\F(V),M)$ directly in order to determine the effect of the map $E(1-P^0)$ in the short exact sequence.  The computation shows:
\begin{align*}
D^0_{U\Bbar/A}(E\F(V),M) & = \Der_{U\Bbar/A}(E\F(V),M)\\
& = \Alg_{U\Bbar \downarrow A}(E\F(V), A\ltimes M)\\
& \cong gr\Mod_{\fpbar \downarrow A}(V, A\ltimes M)\\
& \cong gr\Mod_{\fpbar \downarrow A}(\fpbar\otimes_{\fp}V_0, A \ltimes (\fpbar\otimes_{\fp}M_0) )\\
& \cong gr\Mod_{\fp \downarrow A_0}(V_0, \fpbar \otimes_{\fp} (A_0 \ltimes M_0) ).
\end{align*}
The endomorphism in the short exact sequence (\ref{4term}) is induced by the map  $1-P_0$, which zero except on $\fpbar$, where $P^0$ acts as the frobenius.  Note that $1-P^0$ is surjective on $\fpbar$ as it is the generator of the Galois group of $\fpbar$ over $\fp$.  Thus the kernel of the map on $\fpbar$ is $\fp$, and the cokernel is zero.

The map being surjective shows that \[D^{1}_{U\Bbar/\Fbar\F_0(V)}(\Fbar\F_0(V),M)= \mathrm{Coker} (E(1-P^0)) =0.\]  The kernel being $\fp$ gives rise to the following identifications:
\begin{align*}
D^{0}_{U\Bbar/\Fbar\F_0(V)}(E\F_0(V),M) & \cong \mathrm{Ker}(E(1-P^0))\\
 & \cong gr\Mod_{\fp \downarrow A_0}(V_0, A_0 \ltimes M_0)\\
& \cong \Alg_{U\A \downarrow A_0}(E\F_0V_0, A \ltimes M_0)\\
& \cong \Der_{U\A/A_0}(E\F_0V_0,M_0),
\end{align*} which is the second point.
\end{proof}

\section{The unstable Adams spectral sequence \label{ua}}
The unstable $\hfp$ based Adams spectral sequence abuts to the unstable homotopy groups of a space localized with respect to the ring spectrum $\hfp$.  The following description comes out \cite{bk-uass}, which is described in terms of simplicial sets and is sufficient for our purposes here.  For analogous constructions defined for more general cohomology theories and on the category of spaces rather than simplicial sets, see \cite{etal}.  These should be considered the references for details on the construction and convergence of the unstable Adams spectral sequence.

Fix a simplicial set $X$.  Note that there is an adjunction between simplicial sets and simplicial $\fp$-algebras:
\begin{align*}
\xymatrix{F: \T \ar@<.5ex>[r] &  \ar@<.5ex>[l] s\Alg_{\fp} : U.}
\end{align*}
The right adjoint $U$ is the forgetful functor.  The left adjoint $F$ is the free functor sending $X\mapsto \fp(X)$ to the free $\fp$-algebra on the simplicial set $X$.  The composition $U\circ F$ is a monad on $\T$ and we are interested in the cosimplicial space constructed with this monad.  There are two coface maps \[\xymatrix{\fp(X) \ar@<.5ex>[r] \ar@<-.5ex>[r] & \fp(\fp(X))}\] given by $F(\fp(-))$ and $\fp(F(-))$ respectively.  There is one codegeneracy map \[\xymatrix{\fp(X) &\ar[l] \fp(\fp(X))}\] given by the algebra multiplication of $\fp$ with itself.  Thus iterating this process, we obtain a cosimplicial space that we denote $\fp^{\bullet +1}(X)$ where $\bullet = 0,1,\dotsc$, and $\fp^m(X)$ denotes the $m$ fold composition of the monad $U\circ F$.  In particular, viewing $X$ as the constant cosimplicial object, we obtain a cosimplicial map \begin{align}\label{cosimplicial}X \rightarrow \fp^{\bullet +1}(X)\end{align} whose totalization is the Bousfield--Kan $p$-completion.  When $X$ is nilpotent, as we suppose, this is equivalent to the $\hfp$-localization.  Thus if we have a mapping space $\T(Y,X)$ of nilpotent spaces, we obtain a cosimplicial space $\T(Y,\fp^{\bullet +1}(X))$ whose totalization is $\T(Y,X_{\hfp})$.

\begin{theorem}\label{uass}
There is a Bousfield--Kan spectral sequence abutting to the unstable homotopy groups $\pi_*\T(Y, X_{\hfp})$ whose $E_2$ term is given by:
\[ E_2^{s,t}= \pi^s \pi_t \T(Y,\fp^{\bullet +1}(X)).\]
\end{theorem}
Convergence of this fringed spectral sequence is a delicate issue.  For proof of this theorem and conditions for convergence, see \cite{bk-uass}.  It is interesting to note that to guarantee convergence of this spectral sequence, one need only suppose a space is $p$-complete and nilpotent \cite{pprofinite}.

\subsection{Identification of $E_2$ term}
We are interested in identifying the $E_2$ term in the case that $X$ and $Y$ are of finite type.  The reason we wish to suppose this is that it allows us to describe the $E_2$ term of the spectral sequence in terms of the cohomology of $X$ and $Y$ rather than the homology.  Throughout this section, we suppose that $X$ and $Y$ are of finite type.  These conditions are also necessary for the proof of the main theorem, thus nothing is lost in supposing them here.

Recall the Steenrod algebra as $\A$.  Let $\Alg_{U\A}$ denote the category of unstable $\A$-algebras.  There is an adjoint pair of functors:
\begin{align}\xymatrix{
\F_0: gr\Mod_{\fp} \ar@<.5ex>[r] & \ar@<.5ex>[l] U\A: U
}\end{align}  Here $\F_0$ is the free unstable algebra functor described in \cite{schwartz}.  

We can identify $\pi^*\fp(X) \cong \hfp^*X$, and furthermore 
\begin{lemma}
When $X$ is of finite type, there is an isomorphism of unstable $\A$-algebras:
\begin{align*} \hfp^*\fp(X) \cong \F_0 (\hfp^*X) .\end{align*}
\end{lemma}  Proof follows from a colimit argument applied to the case that $X$ is a sphere.  In this case the lemma says that the cohomology $\hfp^*K(\fp,n)$ is a the free unstable $\A$-algebra on a single generator in degree $n$.

Let $f \in \pi_0\T(Y,X)$ be a basepoint for $\T(Y,X)$.  Then consider 
\begin{align*}
\pi_t\T(Y,X) & \cong \pi_0\T(S^t \wedge Y,X).
\end{align*}  
There is a Hurewitcz-type homomorphism to the category of unstable algebras over the Steenrod algebra $\A$:
\begin{align}\label{hur}
\pi_0\T(S^t \wedge Y,X) \rightarrow & \Hom_{U\A}(\hfp^*X,\hfp^*(S^t \wedge Y)).
\end{align}
Note that the map (\ref{hur}) is an isomorphism if $X$ is of the form $\fp(X)$.  Thus we have an identification of the $E_2$ term of our spectral sequence:
\begin{align*}
\pi^s \pi_t \T(Y,\fp^{\bullet +1}) & \cong \pi^s \pi_0 \T(S^t\wedge Y, \fp^{\bullet +1}(X))\\
& \cong \pi^s \Hom_{U\A}(\hfp^*\fp^{\bullet +1}(X), \hfp^*S^t \wedge Y).
\end{align*}
Since $\hfp^*\fp^{\bullet +1}(X) \cong \F_0 \hfp^*\fp^{\bullet}(X)$ is a simplicial resolution of $\hfp^*X$ by free $U\A$-algebras, we can identify $\pi^s \Hom_{U\A}(\hfp^*\fp^{\bullet +1}(X), \hfp^*S^t \wedge Y)$ as 
\begin{align*} \Ext^s_{U\A}(\hfp^*X,\hfp^*\Sigma^tY) & \qquad t>0, s\geq 1\\
\Hom_{U\A}(\hfp^*X,\hfp^*Y) &\qquad (s,t)=(0,0).
\end{align*}  See \cite{g-aqbk} for a discussion of this identification.  When $t>0$, we can further identify
\begin{align*}
\Ext^s_{U\A}(\hfp^*X,\hfp^*\Sigma^tY) = D^s_{U\A/\hfp^*Y}(\hfp^*X, \Omega^t\hfp^*Y)
\end{align*}
 where $D^*_{U\A/\hfp^*Y}(-,-)$ is an Andr\'e--Quillen type cohomology theory denoting the derived functors of derivations of unstable $\A$-algebras augmented over $\hfp^*Y$.  Note that we are using the fact that Steenrod operations act trivially on the cohomology of suspension spaces.  Observe that $\hfp^*X$ maps to $\hfp^*Y$ via the basepoint map $\hfp^*f$.  Note that $\Omega^t$ refers to a cohomological degree shift.  

\section{The Goerss--Hopkins spectral sequence \label{gh}}
The Goerss--Hopkins spectral sequence is a spectral sequence that computes the homotopy groups of a space of maps between $E_{\infty}$-algebras.  In its most general form, we fix an $\s$-algebra $A$.  Let $E$ be an $A$-algebra such that $\pi_*(E\wedge_A E)$ is flat over $\pi_*E$.  Given $A$-algebras $M, N$ and a map of algebras $\phi: M \rightarrow N$ giving $N$ the structure of an $M$-module, the Goerss--Hopkins spectral sequence computes \[\pi_*(\Alg_A(M,N^{\wedge}_E),\phi)\] in terms of $E_*M$ and $E_*N$.  For details on the construction, see \cite{gh,ghlong,gh-aq}.  We are interested in the case where both $A$ and $E$ are the Eilenberg--Maclane spectrum $\hfpbar$.  In this section, we describe the $E_2$ term of this spectral sequence and identify it algebraically as the left derived functors of derivations of simplicial algebras over an operad.

Let $\Mod_{\hfpbar}$ denote the model category of $\hfpbar$-modules.  Let $\op$ denote the $E_{\infty}$ commutative operad in spaces: each space $\op(m)$ is a contractible, free $\Sigma_m$-space.  Define $\Alg_{\hfpbar} := \Alg_{\op}(\Mod_{\hfpbar})$ to be the category of algebras over the commutative operad in $\hfpbar$-modules.  There is a Quillen adjunction: \[
\xymatrix{\Mod_{\hfpbar} \ar@<.5ex>[r]^{\op(-)} & \Alg_{\hfpbar} \ar@<.5ex>[l],}
\] with right adjoint the forgetful functor and left adjoint $\op(-)$ the free commutative algebra functor.   
The model category structure on $\Alg_{\hfpbar}$ is induced by the model category structure on $\Alg_{\s}$ by identifying $\Alg_{\hfpbar}$ with the under category $\hfpbar \downarrow \Alg_{\s}.$  The model category structure on $\Alg_{\s}$ is induced by the model category structure on a good model for spectra.  One may take the model category structure on spectra of $\s$-modules in \cite{ekmm}, or one can take the positive projective model category structure on symmetric spectra \cite{hss}.  The model structure on $\Alg_{\hfpbar}$ is determined by the fibrations and weak equivalences, which are the maps that are fibrations and weak equivalences on the underlying category of spectra.

In the homotopy category Ho$\Mod_{\hfpbar}$, let $\p$ denote the class of projectives given by the finite cell $\hfpbar$-modules.  This is the smallest class of spectra has the following properties:
\begin{itemize}
\item $\hfpbar \in \p$
\item $\p$ is closed under suspension and desuspension
\item $\p$ is closed under finite wedges
\item For all $P \in \p$, and all $X \in \Mod_{\hfpbar}$, the K\"unneth map \[[P,X]
\rightarrow \Mod_{\fpbar}(\pi_*P,\pi_*X)\] is an isomorphism.
\end{itemize}
Note that this is equivalent to the class of finite coproducts of suspensions and desuspensions of $\hfpbar$.
This allows us to define the $\p$-projective model category structure on $\sMod_{\hfpbar}$ \cite{bousinj,ghlong}.  We say that a map $X_{\bullet} \rightarrow Y_{\bullet}$ in $\sMod_{\hfpbar}$ is  
\begin{itemize}
\item a $\p$-weak equivalence if $[P,X_{\bullet}] \rightarrow [P,Y_{\bullet}]$ is a weak equivalence of simplicial $\fpbar$-modules for all $P\in \p$.
\item a $\p$-fibration if it is a Reedy fibration and $[P,X_{\bullet}] \rightarrow [P,Y_{\bullet}]$ is surjective for all $P\in \p$.
\item a $\p$-cofibration if it has the left lifting property with respect to all $\p$-fibrations that are also $\p$-weak equivalences.
\end{itemize}
This model structure is cofibrantly generated; see \cite{ghlong}.  
Since the $\p$-projective model structure on $s\Mod_{\hfpbar}$ is cofibrantly generated, we obtain a cofibrantly generated model structure on $\sAlg_{\hfpbar}$, which we also call the $\p$-projective model category structure, whose weak equivalences and fibrations are those that are $\p$-weak equivalences and $\p$-fibrations on the underlying category $\sMod_{\hfpbar}$ \cite{hss}.  

\begin{lemma} [{\cite[Proposition 1.4.11]{ghlong}}] \label{funcres}
There exist functorial cofibrant replacements $p: P(X_{\bullet}) \tilde{\twoheadrightarrow} X_{\bullet}$ in the category of simplicial (commutative) $\hfpbar$-algebras such that the degeneracy diagram of $P(X_{\bullet})$ has the form $\op(Z_{\bullet})$, where $Z_{\bullet}$ is free as a degeneracy diagram, and each $Z_n$ is a wedge of finite cell $\hfpbar$-module spectra.
\end{lemma}

We can now describe the Goerss--Hopkins spectral sequence.  
\begin{theorem}[{\cite[Theorem 4.3]{gh}}]
Given $A$ and $B$ in $\Alg_{\hfpbar}$, and a map $\phi: A\rightarrow B$, there is a totalization spectral sequence with $E_2$ term 
\[E_2^{s,t} := \pi^s\pi_t (\Alg_{\hfpbar}(P(A), B), \phi \circ p)\Longrightarrow \pi_{t-s}(\Alg_{\hfpbar}(A,B),\phi).\] 
\end{theorem}

\begin{proof}
Given the standard converges hypotheses, this spectral sequence converges to \[\pi_{t-s}\Tot(\sAlg_{\hfpbar}(P(A),B)).\]  The first thing to check is the equivalence of this space to the space we are interested in:
\begin{align*}
\Tot(\Alg_{\hfpbar}(P(A),B) &\simeq \textrm{Nat}(\Delta^{\bullet}, \sAlg_{\hfpbar}(P(A),B))\\
& \simeq \sAlg_{\hfpbar}(P(A),B^{\Delta^{\bullet}})\\
&  \simeq \Alg_{\hfpbar}(|P(A)|,B)\\
& \simeq \Alg_{\hfpbar}(A,B).
\end{align*}  For proof of the standard technical hypotheses guaranteeing the converges of this spectral sequence, see \cite{gh}.  
\end{proof}

\subsection{Identification of the $E_2$ term}
Let  $A$ and $B$ be $\hfpbar$-algebras.  Fix an algebra map $\varphi: A\rightarrow B$.  Note that the $E_2$ term of the spectral sequence is given by $\pi^s\pi_t\Alg_{\hfpbar}(P(A),B).$  We begin by analyzing $\pi_t\Alg_{\hfpbar}(P(A),B).$  There is a Hurewicz-type homomorphism:
\begin{align*}
\pi_t(\Alg_{\hfpbar}(A,B);\varphi) & \cong \pi_0\Alg_{\hfpbar}(A,B^{S^t})\\
& \rightarrow \Alg_{U\Bbar/B_*}(A_*,B_*[x_{-t}]/(x^r_{-t}))\\
& \cong \Der_{U\Bbar}(A_*,\Omega^t B_*).
\end{align*}
We show this map is an isomorphism if $A = \op(Z),$ which will give us the desired identification when we replace $A$ by $P(A)\cong \op(Z_{\bullet})$.  

Let $A= \op(Z)$.  Fix an $\hfpbar$-algebra map $\varphi: \op(Z) \rightarrow B$ as the basepoint.  There is an induced $\hfpbar$-module map $\varphi': Z \rightarrow \op(Z) \rightarrow B$.  

We begin by analyzing
\begin{align*}
\pi_t \Alg_{\hfpbar}(\op(Z), B) \cong \pi_t \Mod_{\hfpbar}(Z, B)
\end{align*} where the basepoint for the space of $\hfpbar-$module maps is given by $\varphi'$.  The homotopy groups of the space of module maps can be thought of as 
\begin{align*} \pi_0{\T_+}(S^t, \Mod_{\hfpbar}(Z, B)).
\end{align*}

Observe that while the category $\Alg_{\hfpbar}$ is not naturally pointed, the category $\Mod_{\hfpbar}$ is.  Thus the mapping space $\Mod_{\hfpbar}(Z,B)$ has a natural basepoint $0$ given by the zero map.  Thus, while we view the pointed space as an object \[\Map_{\hfpbar}(Z,B) \in * \downarrow \T\] in the category of pointed spaces, it is naturally an object in the category $S^0 \downarrow \T$ of spaces under $S^0$.  Here we think of $S^0$ as two points 0 and 1.  The point 0 maps to the zero map $0$, while the point 1 maps to our original point $\varphi '$.  Since forgetting a basepoint is right adjoint to adding a disjoint basepoint, we have an equivalence:
\begin{align*}
\pi_0 (*\downarrow\T)(S^t,\Mod_{\hfpbar}(Z, B)) \cong \pi_0 (S^0 \downarrow \T)(S^t_+, \Mod_{\hfpbar}(Z, B)).
\end{align*}
The category $\Mod_{\hfpbar}$ is tensored and cotensored over $\T$, thus we have the following equivalences:
\begin{align*}
\pi_0 (S^0 \downarrow \T)(S^t_+,\Mod_{\hfpbar}(Z, B)) &\cong \pi_0 \Mod_{S^0 \otimes Z \downarrow \hfpbar}(S^t_+\otimes Z, B)\\
& \cong \pi_0 \Mod_{\hfpbar \downarrow B^{S^0}} (Z, B^{S^t_+})
\end{align*}
Here $Z$ is an $\hfpbar-$module over $B$ via the map $\varphi '$, while $B^{S^t_+}$ is via the map of spaces $S^0 \rightarrow S^t_+$ sending 0 to $+$ and 1 to the natural basepoint 1 in $S^t$.  Observe that since $B$ is actually an $\hfpbar-$algebra, so is $B^{S^t_+}$.  Now we can identify this term as 
\begin{align*}
\pi_0 \Mod_{\hfpbar \downarrow B^{S^0}} (Z, B^{S^t_+}) \cong \Mod_{\fpbar \downarrow B_*} (\pi_*Z, \pi_*B^{S^t_+}).
\end{align*} 

Since $B^{S^t_+}$ is naturally an $\hfpbar-$algebra, its homotopy is naturally the unstable $\Bbar-$algebra $B^*(S^t) \cong B_*[x_{-t}]/(x_{-t}^2)$.  The algebra structure is given by the product in the cohomology theory $B^*$.  Let $(n,m)$ and $(n',m')$ be elements in $B_* \oplus \Omega^tB_*$, then \[(n,m)\cdot(n',m') = (nn',nm'+n'm)\] since $mm' = 0$ as there are no higher cohomology groups.  Thus we might more appropriately write this algebra as the square-zero extension  $B^* \ltimes \Omega^tB^*$.  Using the left adjoint to the forgetful functor from $U\Bbar$ to $\Mod_{\fpbar}$ we have an identification
\begin{align*}
\Mod_{\fpbar \downarrow B_*} (\pi_*Z, B_*\ltimes \Omega^tB_*) \cong \Alg_{U\Bbar \downarrow B^*} (E\F(\pi_*Z), N^*\ltimes \Omega^tN^*), 
\end{align*}
proving the map in the Hurewicz type homomorphism is in fact an isomorphism.

By Lemma \ref{adaption}, $\pi_*P(A) \cong \pi_*\op(Z) \cong \Fbar\F(\pi_*Z) .$  By the previous remark, we see that:
\begin{align}\label{huriso}
\pi_t \Alg_{\hfpbar}(P(A), B) & \cong \pi_t \Alg_{\hfpbar}(\op(Z_{\bullet}), B) \cong \Der_{U\Bbar}(\pi_*\op(Z_{\bullet}),\Omega^tB_*).
\end{align}
Observe that $\pi_*\op(Z_{\bullet}) \rightarrow A_*$ is a cofibrant simplicial resolution of $A_*$ by levelwise free unstable $\Bbar$-algebras.  Taking $\pi^s$ of the cosimplicial complex (\ref{huriso}) we obtain the desired identification of the $E_2$ term as the right derived functors of derivations of unstable $\Bbar$ algebras.

\begin{proposition}
The $E_2$ term of the Goerss--Hopkins spectral sequence can be identified as
\begin{align*}
\pi^s\pi_t \Alg_{\hfpbar}(A,B)  & \cong D^s_{U\Bbar}(A_*,\Omega^tB_*) & \qquad t>0\\
& \cong \Alg_{U\Bbar}(\hfpbar^*X,\hfpbar^*Y) & \qquad s = t = 0.
\end{align*}
\end{proposition}

Given a pointed, p-complete, nilpotent space $X$, and a nilpotent space $Y$, we are interested in computing the space $\Alg_{\hfpbar}(\hfpbar^{\Sigma^{\infty}X_+},\hfpbar^{\sinf Y_+})$.  The of the description of the Goerss--Hopkins spectral sequence computing the homotopy groups of this spaces simplifies.  
\begin{proposition} \label{ghss}  Given spaces $X$ and $Y$, the Goerss--Hopkins spectral sequence computing the homotopy groups of $\psi(X)$ has $E_2$ term given by:
\begin{align*}E_2^{s,t}= &D^s_{U\Bbar/\hfpbar^*Y}(\hfpbar^* X, \Omega^t\hfpbar^*Y)  & t>0\\
E_2^{0,0} = & \Alg_{U\Bbar}(\hfpbar^*X,\hfpbar^*Y) &
\end{align*}
abutting to
\begin{align*}
\pi_{t-s}( \Alg_{\hfpbar}(\hfpbar^{\sinf X_+}, \hfpbar^{\sinf Y_+})).\end{align*}
\end{proposition}

\section{Comparing cosimplicial spaces \label{alg}}

Recall that we have the map $\psi$ from (\ref{mainmap}).  In particular, $\psi$ induces a map of cosimplicial spaces:
\begin{align}\label{maineqn}
\T(Y,\fp^{\bullet+1}(X)) \rightarrow \Alg_{\hfpbar}(\hfpbar^{\Sigma^{\infty}\fp^{\bullet+1}(X)_+}, \hfpbar^{\sinf Y_+})).
\end{align} 

The goal is to prove that these two cosimplicial spaces are weakly equivalent.  Thus after applying $\pi^s\pi_t$ to this map, we will obtain an isomorphism between the $E_2$ term of the unstable Adams spectral sequence and the Goerss--Hopkins spectral sequence.

\begin{lemma}\label{mainres}
The map 
\begin{align}\label{thesimpaug}
\hfpbar^{\Sigma^{\infty}\fp^{\bullet+1}(X)_+} \rightarrow \hfpbar^{\Sigma^{\infty}X_+}
\end{align} is a $\p$-weak equivalence in $\sAlg_{\hfpbar}$.   \end{lemma}
\begin{proof}
Recall that \begin{align*}
\hfpbar^{\Sigma^{\infty}\fp^{\bullet+1}(X)_+} \rightarrow \hfpbar^{\Sigma^{\infty}X_+}
\end{align*} is a $\p$-weak equivalence if 
\begin{align*} \big[P,\hfpbar^{\Sigma^{\infty}\fp^{\bullet+1}(X)_+} \big] \rightarrow \big[P,\hfpbar^{\Sigma^{\infty}X_+}\big]
\end{align*} is a weak equivalence for all $P = \op(Z)$, where $Z$ is a finite cell $\hfpbar$-module.  Since this is the same as saying 
\begin{align*} \Mod_{\fpbar}(\pi_*P,\pi_*\hfpbar^{\Sigma^{\infty}\fp^{\bullet+1}(X)_+}) \rightarrow \Mod_{\fpbar}(\pi_*P,\pi_*\hfpbar^{\Sigma^{\infty}X_+}) \end{align*} is a weak equivalence of simplicial abelian groups for all $P$ in $\p$, and moreover $\pi_*P$ is projective as an $\fpbar$-module, it suffices to check that
\[\xymatrix{ \hfpbar^*\fp^{\bullet+1}(X) \ar[r] &\hfpbar^*X
}\] is a weak equivalence of simplicial abelian groups.  But asking whether $X\rightarrow \fp^{\bullet+1}(X)$ is an $\hfp$ equivalence is asking whether $\fp(X) \rightarrow \fp^{\bullet+2}(X)$ is a weak equivalence. This follows since there is an extra degeneracy giving rise to a contracting homotopy; see for example \cite{bk-uass,gj}.  \end{proof}

Write $W_{\bullet} \rightarrow W$ for the simplicial resolution (\ref{thesimpaug}).  For each $n$ we obtain a functorial simplicial $\p$-resolution $P(W_n)\rightarrow W_n.$  Thus we obtain a bisimplicial object $P(W_{\bullet})$.  

There is an injective model category structure on bisimplicial $\hfpbar$-algebras {\cite[Proposition A.2.8.2] {topos}} so that a map $A_{\bullet \bullet} \rightarrow B_{\bullet \bullet}$ is:
\begin{itemize}
\item a cofibration if $A_{n \bullet} \rightarrow B_{n \bullet}$ is a $\p$-cofibration for all $n$,
\item a weak equivalence if $A_{n \bullet} \rightarrow B_{n \bullet}$ is a $\p$-equivalence for all $n$, 
\item and fibrations are determined.
\end{itemize}  Call this model category structure the injective$-\p$ model structure.  The idea is that in the vertical simplicial direction, the cofibrations and weak equivalences are determined levelwise. 
\begin{lemma}
In the injective$-\p$ model structure, $P(W_{\bullet})$ is cofibrant, and  $P(W_{\bullet}) \rightarrow W_{\bullet}$ is a weak equivalence, where $W_{\bullet}$ is viewed as a vertically constant bisimplicial $\hfpbar$-algebra.  
\end{lemma}

\begin{proof}
This follows since we need only check that levelwise that $P(W_n)$ is $\p$-cofibrant and $P(W_n) \rightarrow W_n$ is a $\p$-weak equivalence, which is true by construction.  
\end{proof}  Define $\widetilde{W_{\bullet}} := \textrm{diag}(P(W_{\bullet}))$.  

\begin{lemma}\label{cof}
The diagonal object $\widetilde{W_{\bullet}}$ is $\p$-cofibrant, and $\widetilde{W_{\bullet}} \rightarrow W_{\bullet}$ is a $\p$-weak equivalence.  
\end{lemma}

\begin{proof}
The proof is obtained by dualizing the proof found in {\cite[Lemmas 6.9 and 6.10] {bousinj}}.
\end{proof}

We can compose the map (\ref{maineqn}) with the edge map 
\[ \Alg_{\hfpbar}(W_{\bullet}, \hfpbar^{\sinf Y_+})) \rightarrow \Alg_{\hfpbar}(P(W_{\bullet}), \hfpbar^{\sinf Y_+}))\] and we obtain a diagram of cosimplicial spaces:
\begin{align}\label{diag} \xymatrix{
\T(Y,\fp^{\bullet+1}(X)) \ar[r] & \mathrm{diag}\Alg_{\hfpbar}(P(W_{\bullet}), \hfpbar^{\sinf Y_+}))\\
& \Alg_{\hfpbar}(P(\hfpbar^{\sinf X_+}), \hfpbar^{\sinf Y_+})). \ar[u]
}\end{align}  The goal is to prove that when $X$ and $Y$ are of finite type, these maps are all equivalences.

Before we get to the main result of this section, we take a moment to resolve homotopy groups of the spectra $W_{n}$ as a colimit of free unstable $\B$-algebras over $\fpbar$.  Define \[W_{-1}:= \hfpbar^{\sinf X_+}.\]  The proof of Proposition \ref{coeq} applied to the spectra in our simplicial resolution gives rise to cofiber sequences:
\[\xymatrix{E\F(\pi_*W_{n-1}) \ar[r]^{1-P^0}  &  E\F(\pi_*W_{n-1}) \ar[r] 
&E\F_0(\pi_*W_{n-1}) \cong \pi_*W_n} \]
\noindent for each $n \geq 0$.  

\begin{proposition}\label{totss}
The totalization spectral sequence
\begin{align*}\pi^s\pi^u\pi_t\sAlg_{\hfpbar}(P(W_{\bullet}),\hfpbar^{\sinf Y_+}) \Longrightarrow \pi^{u+s}\pi_t\sAlg_{\hfpbar}(\widetilde{W_{\bullet}},\hfpbar^{\sinf Y_+}), \end{align*} collapses, and induces an isomorphism
\begin{align*}
\pi^s\pi^0\pi_t\sAlg_{\hfpbar}(P(W_{\bullet}),\hfpbar^{\sinf Y_+}) \rightarrow &\pi^{s}\pi_t\sAlg_{\hfpbar}(\widetilde{W_{\bullet}},\hfpbar^{\sinf Y_+}) \end{align*}
\end{proposition}

 \begin{proof}
Since $\widetilde{W_{\bullet}}$ is cofibrant, and its realization is weakly equivalent $W$, the right hand side is the $E_2$ term of the Goerss--Hopkins spectral sequence in Proposition \ref{ghss}.  Observe that for each $n$ and $t>0$, the left hand side can be identified as
\begin{align*}
\pi^s\pi^u\pi_t\sAlg_{\hfpbar}(P(W_n),\hfpbar^{\sinf Y_+}) =\pi^s D^u_{U\B/\hfpbar^*Y}(\pi_*W_n, \Omega^t\hfpbar^{*}Y).
\end{align*}  

Since $\pi_*W_n$ is isomorphic to $E\F_0(\pi_*W_{n-1})$, we can apply Theorem \ref{cor}.  The consequence of this theorem is that:
\begin{align*}
D^u_{U\Bbar/\hfpbar^*Y}(\pi_*W_n, \hfpbar^{*-t}Y) &= 0 
\end{align*} if  $u \neq 0,$  and 
\begin{align*}
D^0_{U\Bbar/\hfpbar^*Y}(\pi_*W_n, \hfpbar^{*-t}Y) = D^0_{U\A/\hfp^*Y}(\pi_*W_{n-1}, \Omega^t\hfp^{*}Y).
\end{align*}
\end{proof}  

\begin{proposition}\label{eq}
The map $k: \widetilde{W_{\bullet}} \rightarrow W_{\bullet}$ induces isomorphisms \[
\xymatrix{\pi^s\pi_t\big( \Alg_{\hfpbar}(W_{\bullet},\hfpbar^{\sinf Y_+});\varphi \big) \ar[r] & \pi^s\pi_t\big(\Alg_{\hfpbar}(\widetilde{W_{\bullet}}, \hfpbar^{\sinf Y_+}); k\circ \varphi \big).}\]
\end{proposition}

\begin{proof}
The proof of Proposition \ref{totss} exhibits an isomorphism 
\[ \xymatrix{\pi^s\pi_t\Alg_{\hfpbar}(W_{\bullet},\hfpbar^{\sinf Y_+}) \ar[r]^-{\simeq} &\pi^s\pi^0\pi_t\sAlg_{\hfpbar}(P(W_{\bullet}),\hfpbar^{\sinf Y_+}). }\]  The collapse of the spectral sequence in Proposition \ref{totss} gives an isomorphism 
\begin{align*} \xymatrix{\pi^s\pi^0\pi_t\Alg_{\hfpbar}(P(W_{\bullet}),\hfpbar^{\sinf Y_+}) \ar[r]^-{\simeq} & \pi^s\pi_t\Alg_{\hfpbar}(\widetilde{W}_{\bullet},\hfpbar^{\sinf Y_+}).} \end{align*}

The composition yields the desired isomorphism.  The fact that it is coming from the map $\widetilde{W}_{\bullet} \rightarrow W$ follows from the fact that the two isomorphisms are induced by maps in the diagram of bisimplicial $\hfpbar$-algebras: 
\[ \xymatrix{ 
P(W_{\bullet})  \ar[dr]_{\simeq} &&\ar[ll]_{\simeq} \widetilde{W}_{\bullet} \ar[dl]\\
&W_{\bullet} .} \]  

In particular,  by the two out of three property, the map $ \widetilde{W}_{\bullet} \rightarrow  W_{\bullet}$ is an injective $\p$-weak equivalence of constance bisimplicial spectra, and thus is a $\p$-weak equivalence of simplicial spectra.  Thus the map is a $\p$-weak equivalence out of a $\p$-cofibrant object.
\end{proof}

\begin{theorem}\label{thepoint}
When $X$ and $Y$ are finite type, nilpotent spaces, the diagram (\ref{diag}) is a diagram of weak equivalence of cosimplicial spaces.  Moreover, these maps induce weak equivalences
\begin{align*}
\xymatrix{
\T(Y,X_p) \ar[r] & \Tot \Alg_{\hfpbar}(\widetilde{W}_{\bullet},\hfpbar^{\sinf Y_+})\\
& \Alg_{\hfpbar}(\hfpbar^{\sinf X_+}, \hfpbar^{\sinf Y_+}). \ar[u]
}\end{align*}
\end{theorem}
\begin{proof}
Consider the map
\[ \pi^s\pi_t\T(Y,\fp^{\bullet+1}(X)) \rightarrow \pi^s\pi_t\Alg_{\hfpbar}(\widetilde{W}_{\bullet}, \hfpbar^{\sinf Y_+}).\]  From Proposition \ref{eq}, we know that the map \[\pi^s\pi_t \Alg_{\hfpbar}(W_{\bullet},\hfpbar^{\sinf Y_+}) \rightarrow \pi^s\pi_t\Alg_{\hfpbar}(\widetilde{W}_{\bullet}, \hfpbar^{\sinf Y_+})\] is an isomorphism.  Furthermore, we know that the map \[\pi^s\pi_t \T(Y,\fp^{\bullet+1}(X)) \rightarrow \pi^s\pi_t \Alg_{\hfpbar}(W_{\bullet},\hfpbar^{\sinf Y_+}) \] is an isomorphism when $X$ and $Y$ are of finite type since both sides can be identified as derived functors of derivations of unstable $\A$-algebras over $\hfp^*Y$.  

Similarly, Proposition \ref{totss} exhibits an isomorphism:
\[\pi^s\pi_t\Alg(\widetilde{W}_{\bullet},\hfpbar^{\sinf Y_+}) \leftarrow \pi^s\pi_t\Alg(P(\hfpbar^{\sinf X_+}),\hfpbar^{\sinf Y_+}).\]

Thus by \cite{bous}, since there is an isomorphism of the $E_2$ terms of these spectral sequences, there is an induced weak equivalence on the total spaces in the diagram (\ref{diag}).  

In particular, by Lemma \ref{cof}, the map $\widetilde{W}_{\bullet} \rightarrow W$ is a $\p$-weak equivalence, and $\widetilde{W}_{\bullet}$ is $\p$-cofibrant.  Given the functorial cofibrant replacement $P(W)$ of $W$ guaranteed by Proposition \ref{funcres}, $P(W) \rightarrow W$ is a trivial $\p$-fibration.  Thus model category theory guarantees a $\p$-weak equivalence $w: \widetilde{W}_{\bullet}\rightarrow P(W)$ between cofibrant objects.  Thus the associated cosimplicial mapping spaces can be naturally identified.
\end{proof}

\noindent \emph{Remark:}  This give an alternate proof of Mandell's theorem (Theorem \ref{mandell}).  However, the key computations are the same as in his original proof.

\begin{corollary}\label{maincor}
When $X$ and $Y$ are finite type, nilpotent spaces, the map \[\T(Y,\fp^{\bullet+1}(X)) \rightarrow \Alg_{\hfpbar}(\hfpbar^{\sinf \fp^{\bullet+1}(X)_+}, \hfpbar^{\sinf Y_+})\] induces and isomorphism between the unstable Adams spectral sequence and the Goerss--Hopkins spectral sequence.
\end{corollary}
\begin{proof}
Theorem \ref{thepoint} exhibits an isomorphism when applying $\pi^s\pi_t$ to both of these cosimplicial spaces.  Since $\pi^s\pi_t\T(Y, \fp^{\bullet+1}(X))$ can be identified as the $E_2$ term of the unstable Adams spectral sequence, all that is left is to identify \[\pi^s\pi_t\ \Alg_{\hfpbar}(\hfpbar^{\sinf \fp^{\bullet+1}(X)_+}, \hfpbar^{\sinf Y_+})\] as the $E_2$-term of the Goerss--Hopkins spectral sequence.  

Proposition \ref{eq} exhibits an isomorphism \[\pi^s\pi_t\ \Alg_{\hfpbar}(\hfpbar^{\sinf \fp^{\bullet+1}(X)_+}, \hfpbar^{\sinf Y_+}) \cong \pi^s\pi_t\ \Alg_{\hfpbar}(\widetilde{W}_{\bullet}, \hfpbar^{\sinf Y_+}).\] 

By Lemma \ref{cof}, the map $\widetilde{W}_{\bullet} \rightarrow W$ is a $\p$-weak equivalence, and $\widetilde{W}_{\bullet}$ is $\p$-cofibrant.  Given the functorial cofibrant replacement $P(W)$ of $W$ guaranteed by Proposition \ref{funcres}, $P(W) \rightarrow W$ is a trivial $\p$-fibration.  Thus model category theory guarantees a $\p$-weak equivalence $w: \widetilde{W}_{\bullet}\rightarrow P(W)$ between cofibrant objects.  Thus the associated cosimplicial mapping spaces can be naturally identified.  In particular, applying $\pi^s\pi_t$ of either space can be naturally identified as the $E_2$ term of the Goerss--Hopkins spectral sequence.

Thus Proposition \ref{eq} truly exhibits an isomorphism between the unstable Adams spectral sequence and the Goerss--Hopkins spectral sequence.
\end{proof}

%\bibliography{biblio}
%\bibliographystyle{plain}
\end{document}